%% file: CriticalMajorant.tex
\documentclass[a4paper,12pt]{amsart}
\usepackage{standalone}
\usepackage{amsfonts}
\usepackage{amsthm}
\usepackage{amssymb}
\usepackage{amsmath}
\usepackage{theoremref}
\usepackage{enumerate}
\usepackage{bbm}
\usepackage{bm}
\usepackage{pgfplots}
\pgfplotsset{compat=1.18} 
\usepgfplotslibrary{fillbetween}
\usetikzlibrary{intersections}
\usepackage{mathtools}
\usepackage{tkz-euclide}
\usetikzlibrary{patterns}
\usepackage{hyperref}

\numberwithin{equation}{section}

\newcommand{\C}{\mathcal{C}}

\newcommand{\T}{\mathbb{T}}

\newcommand{\conj}[1]{\overline{#1}}
\newcommand{\D}{\mathbb{D}}

\newcommand{\R}{\mathbb{R}}
\renewcommand{\H}{\mathbb{H}}
\newcommand{\Po}{\mathcal{P}}

\newcommand{\hb}{\mathcal{H}(b)}

\renewcommand\Re{\operatorname{Re}}
\renewcommand\Im{\operatorname{Im}}

\newtheorem{thm}{Theorem}[section]
\newtheorem*{thm*}{Theorem}
\newtheorem{lem}[thm]{Lemma}

\newtheorem*{cor*}{Corollary}
\newtheorem{prop}[thm]{Proposition}
\theoremstyle{definition}

\theoremstyle{definition}

\newtheorem{defn}[thm]{Definition}
\newtheorem*{defn*}{Definition}

\newtheorem{claim*}{Claim}

\title[A critical majorant for the Khinchin-Ostrowski property]{A critical majorant for the Khinchin-Ostrowski property}

\author{Bartosz Malman}

\address{Division of Mathematics and Physics, 
        Mälardalen University,
		721 23 Västerås, Sweden}
\email{bartosz.malman@mdu.se} 

\begin{document}

\begin{abstract}
    In this short note we prove an optimal version of a classical result. Given a majorant determining a growth restriction on functions in the unit disk $\D$, we say that a set $E$ on the unit circle $\T$ is a \textit{uniqueness set}, or has the \textit{Khinchin-Ostrowski property}, with respect to the majorant, if any sequence of analytic polynomials satisfying the growth restriction which converges in an appropriate sense to $0$ on $E$, in fact is forced to converge to $0$ in $\D$ also. Theorems proved by Kegejan and Khrushchev state that if $E$ has positive Lebesgue measure and satisfies a generalized Beurling-Carleson condition, then for an appropriate majorant the Khinchin-Ostrowski property is satisfied. A technical point in Khrushchev's proof is the estimation of the harmonic measure in a Privalov-type domain which requires logarithmic integrability of the majorant. This forbids the application of his result to certain types of generalized Beurling-Carleson conditions. Here, we dispose of the integrability assumption on the majorant. To do so, we use a Joukowski-Privalov domain which is obtained by removing from the unit disk the areas enclosed by hyperbolic geodesics between the endpoints of intervals complementary to $E$. For this type of domain the method of Khrushchev applies, but the harmonic measure may be estimated more accurately by simple explicit formulas for conformal mappings. As a consequence, we find the critical majorant at which Beurling-Carleson type conditions stop determining the Khinchin-Ostrowski property of a set, and above which the containment of intervals is the only relevant characteristic. We discuss also weighted versions of the Khinchin-Ostrowski property, and apply our result to establish the remarkable precision of a one-sided spectral decay condition which detects the local logarithmic integrability of a function.
\end{abstract}

\thanks{Part of this work has been done during a visit of the author at the Department of Mathematics and Statistics, University of Reading, funded by the University of Reading and the London Mathematical Society. He gratefully acknowledges the financial support. The author is also partially supported by Vetenskapsrådet (VR2024-03959).}

\maketitle

\section{Introduction}

\subsection{Main result} 
\label{S:MainResultSec}
Let $E$ be a closed subset of the unit circle $\T := \{ z \in \mathbb{C} : |z| = 1\}$. If $(p_n)_n$ is a sequence of analytic polynomials in a single complex variable, and $p_n \to 0$ on $E$ in some appropriate sense, then by analyticity of the polynomials involved we may expect for the sequence $(p_n)_n$ to collapse to $0$ also in the unit disk $\D := \{ z \in \mathbb{C} : |z| < 1 \}$. In the literature, one often refers to such results as \textit{Khinchin-Ostrowski properties}, after the classical result on pointwise limits of sequences in the Nevanlinna class of the disk (see \cite[p. 278]{havinbook}). To prove a rigorous statement of this type, one needs a growth restriction on $(p_n)_n$ in the unit disk, and $E$ needs to be sufficiently \textit{massive} to make the null-convergence of $(p_n)_n$ on $E$ influence the sequence sufficiently strongly inside the disk. The purpose of this note is to reconcile the related results of Khrushchev from \cite{khrushchev1978problem} and the author from \cite{malman2023revisiting}, and force the two works to meet in a remarkably sharp statement of the kind just described. 

The article \cite{malman2023revisiting} can be interpreted as a study of the Khinchin-Ostrowski property, and related matters, in the more general context of weights on $\T$ instead of its subsets $E$. Let us state precisely what we mean by the weighted context, and how it specializes to the classical setting. 

\begin{defn}[\textbf{Khinchin-Ostrowski property}]
\thlabel{D:KOprop} Let $\lambda(x)$ be an increasing positive function (a majorant) defined for $x \in [1, \infty)$, and let $w$ be a non-negative Borel function on $\T$. We say that the pair $(\lambda, w)$ satisfies the Khinchin-Ostrowski property, if for every sequence of analytic polynomials $(p_n)_n$ satisfying the two conditions
\begin{enumerate}[(i)]
    \item $|p_n(z)| \leq C \lambda\Big( \frac{1}{1-|z|} \Big), z \in \D$, for some $C > 0$,
    \item $\lim_{n \to \infty} \int_\T |p(z)| w(z) |dz| = 0$,
\end{enumerate}
we also have that $(p_n)_n$ collapses to $0$ inside the disk:\[ \lim_{n \to \infty} p_n(z) = 0, \quad z \in \D.\]
\end{defn}

Setting $w = 1_E$, the indicator function of a set $E$, the above definition reduces the the classical setting (the convergence $p_n \to 0$ on $E$ being interpreted as convergence in $L^1(E)$). In that case, we shall write $(\lambda, E)$ instead of $(\lambda, 1_E)$. The unweighted context is our main concern, but our results are both inspired by and have implications in the weighted context (see the theorems in Section~\ref{S:ApplicationSubsec}).

We will unveil a special role for the Khinchin-Ostrowski property of the following family of exponential majorants ($c > 0$):
\begin{equation} \label{E:ExpGrowthNew}
\lambda\bigg(\frac{1}{1-|z|} \bigg) = \exp \bigg( \frac{c}{1-|z|} \bigg), \quad z \in \D \tag{ExpGrowth}. \end{equation}
A simple consequence of the results in \cite{malman2023revisiting} (and, in the unweighted context $w = 1_E$, of the results of Khrushchev in \cite{khrushchev1978problem}) is that if the integral \[ \int_I \log w(z) \, |dz| \] diverges over every interval $I \subset \T$, and $\lambda$ is as in \eqref{E:ExpGrowthNew}, then the pair $(\lambda, w)$ does \textit{not} satisfy the Khinchin-Ostrowski property. A fortiori, the property is not satisfied if $\lambda$ in \eqref{E:ExpGrowthNew} is replaced by a larger majorant. For instance, under the two conditions in \thref{D:KOprop} we can still ensure that $p_n \to 1$ on $\D$. Reversing the roles of the disk and the circle, in fact $(p_n)_n$ can be made to converge in $L^1(w)$ to any given function in that space, while simultaneously converging to $0$ in $\D$. Note that if $w = 1_E$, then the divergence of the above logarithmic integral is clearly equivalent to $E$ containing no intervals (if $E$ is not closed then we interpret this statement as $|I \setminus E| > 0$ for every interval $I$). So such a set rather spectacularly fails to satisfy the classical Khinchin-Ostrowski property when paired with a majorant in \eqref{E:ExpGrowthNew}. Conversely, if $\log w$ is integrable over some interval (or if $E$ contains an interval), then the Khinchin-Ostrowski property will hold even if the majorant in \eqref{E:ExpGrowthNew} is replaced by a much larger one (see \cite{malman2023revisiting} for details, and much more general weighted results).

Our goal here is to understand exactly for which majorants the Khinchin-Ostrowski property can be described in terms of intervals alone, as in the results mentioned above. According to what was said there, relevant for this question are the majorants which are smaller than the ones in \eqref{E:ExpGrowthNew}. Therefore, let us introduce largest majorants which are strictly smaller than the ones in \eqref{E:ExpGrowthNew}. Consider
\begin{equation}
    \label{E:LessThanExpGrowth}
    \lambda_h\bigg( \frac{1}{1-|z|} \bigg) := \exp \bigg( \frac{h(1-|z|)}{1-|z|} \bigg), \quad z \in \D,
\end{equation} where $h$ is a non-negative continuous function satisfying $h(0) = 0$. In order to avoid trivialities, we assume that $h$ decays slow enough to have \[ \lim_{x \to 0+} h(x)/x = \infty.\] We will also work under the following rather mild regularity assumptions
\begin{equation}
    \label{E:RegCond}
    h(x) \text{ is increasing, and } h(x)/x \text{ is decreasing, in } x \tag{Reg}
\end{equation} but we will impose no further conditions on $h$. In particular, $h(x)$ may tend to $0$ arbitrarily slowly as $x$ decreases to $0$. Note that if $h$ is a non-negative function satisfying $\lim_{x \to 0+} h(x) = h(0) = 0$ but not \eqref{E:RegCond}, then there exists a continuous function $\widetilde{h}$ pointwise larger than $h$ for which \eqref{E:RegCond} is satisfied, and $\widetilde{h}(0) = 0$. For instance, we may take $\widetilde{h}$ to be the least concave majorant of $h$.

Establishing the following result is the aim of this note.
\begin{thm} \thlabel{T:MainTheorem1}
    Let $h$ be any function satisfying the above description. There exists a closed set $E$ of positive Lebesgue measure that contains no intervals and such that the Khinchin-Ostrowski property holds for the pair $(\lambda_h, E)$.
\end{thm}

A more precise statement identifies the critical property of $E$ needed to establish the theorem. An \textit{$h$-Beurling-Carleson set} is a closed subset of $\T$ for which the family of complementary intervals $\C(E)$ satisfies
\begin{equation}
    \label{E:CarlesonCondDef}
    \sum_{\ell \in \C(E)} h(|\ell|) < \infty.
\end{equation} 
Here $|\ell|$ denotes the length of the interval $\ell$. For any $h$ as above, one may readily construct an $h$-Beurling-Carleson set $E$ of positive Lebesgue measure which contain no intervals. We shall see that any such set satisfies the conclusion of \thref{T:MainTheorem1}. 

\subsection{A rationale for this note's existence} We should justify more carefully the reason for why this note was written, for our main result has a very close resemblance to a well-known result of Khrushchev from \cite{khrushchev1978problem}. We explain the difference and what exactly \thref{T:MainTheorem1} accomplishes. 

Khrushchev proved the following theorem (see \cite[Theorem 3.1]{khrushchev1978problem}).

\begin{thm*}[\textbf{Khrushchev's Uniqueness Theorem}]
    \thlabel{T:KhrushchevTheorem}
    Assume that $E$ has positive Lebesgue measure and that the family of complementary intervals $\C(E)$ satisfies
    \begin{equation}
        \label{E:IntegralCarlesonCond}
        \sum_{\ell \in \C(E)} \int_0^{|\ell|} \log \big( \lambda(1/t) \big)\,dt < \infty.
    \end{equation}
    Then $(\lambda, E)$ satisfies the Khinchin-Ostrowski property.
\end{thm*}

Note the implicit integrability hypothesis in \eqref{E:IntegralCarlesonCond}: we must have 
\begin{equation}
    \label{E:IntHypothesis}
    \int_0^1 \log \big(\lambda (1/t)\big) \, dt< \infty
\end{equation} to apply the theorem. Let 
\begin{equation}
    \label{E:hLambdaConnection}
    h(t) = t \log (\lambda (1/t) \big), \quad t \in (0,1],
\end{equation} that is, $\lambda = \lambda_h$ according to \eqref{E:LessThanExpGrowth}. For many naturally appearing functions $h$, such as $h(t) = \text{constant} \cdot t \log(1/t)$ and $h(t) = \text{constant} \cdot t^\alpha$, $\alpha \in (0,1)$, the condition \eqref{E:IntegralCarlesonCond} is equivalent to \eqref{E:CarlesonCondDef} if $h$ and $\lambda$ are connected by the equation \eqref{E:hLambdaConnection}. 

An elementary estimation which uses the integrability \eqref{E:IntHypothesis} and the earlier stated monotonicity condition on $\lambda$ shows that $\lim_{t \to 0+} h(t) = 0$. It follows that the majorants appearing in Khrushchev's result are of the type \eqref{E:LessThanExpGrowth}. But the integrability \eqref{E:IntHypothesis} says that 
\[ \int_0^1 \frac{h(t)}{t} \, dt < \infty\] which is an additional smallness hypothesis of $h$ which prohibits it from decaying too slowly near $0$, and so the full range of majorants in \eqref{E:LessThanExpGrowth} is not covered by those satisfying \eqref{E:IntHypothesis}. One such majorant $\lambda_h$ is corresponding to $h(t) = 1/\log(1/t)$.
The principal upside of our \thref{T:MainTheorem1} is that it is free of any additional integrability or size restrictions on $h$, and thus it applies to \textit{any} majorant in \eqref{E:LessThanExpGrowth}  (but, confessedly, at the cost of the mild regularity conditions in \eqref{E:RegCond}, which may or may not be necessary). \textit{It is precisely this improvement that identifies the family of majorants \eqref{E:ExpGrowthNew} to be critical for the Khinchin-Ostrowski property.} Explicitly, the improved result in \thref{T:MainTheorem1} shows that:

\begin{itemize}
    \item corresponding to every majorant smaller than the ones in \eqref{E:ExpGrowthNew}, there exists a closed set $E$ containing no intervals for which the corresponding Khinchin-Ostrowski property is satisfied,
    \item for any majorant at least as large as the majorants in \eqref{E:ExpGrowthNew}, the set $E$ must contain an interval if the corresponding Khinchin-Ostrowski property is to hold.
\end{itemize}

In the weighted context, the same conclusion holds: the majorants in \eqref{E:ExpGrowthNew} are the critical ones in the context of determining the Khinchin-Ostrowski property for $(\lambda, w)$ from integrability of $\log w$ over intervals alone.

\subsection{The technique: estimation of harmonic measure}
\label{S:TechniqueSection}
In our proof of \thref{T:MainTheorem1}, we will re-use the basic approach of Khrushchev, which itself is similar to one of Kegejan from \cite{kegejanex}. See \cite[Section 3]{khrushchev1978problem} for details of the claims made next. The main difference will be that our principal estimates will be sharper.

\begin{figure}
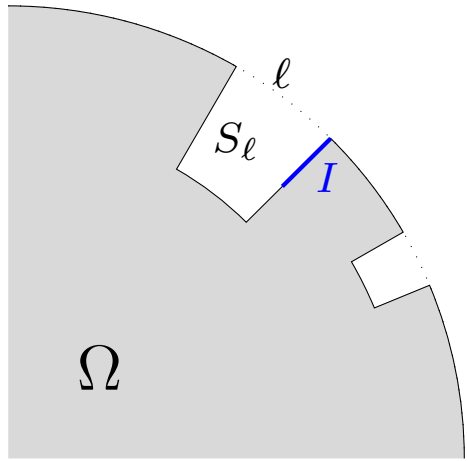

    \centering
    \includestandalone[scale=2]{khrushchevsquare}
    \caption{Domain used in Khrushchev's proof.}
    \label{F:Figure3}
\end{figure}

Let $\omega(\Omega, A, z)$ denote the harmonic measure of a Borel set $A$ contained in $\partial \Omega$ at the point $z$ inside a given domain $\Omega$. The connection of values of $(p_n)_n$ on $E$ with the values in $\D$ is established by constructing a nice simply connected domain $\Omega \subset \D$ which satisfies $\partial \Omega \cap \T = E$. We assume that the boundary $\partial \Omega$ is a rectifiable Jordan curve, and that $0 \in \Omega$. Then the harmonic measure $\omega = \omega(\Omega, \cdot, 0)$ exists, and subharmonicity of $\log |p_n|$ implies readily that 
\begin{equation}
    \label{E:SubharmLogpn}
    \log |p_n(0)| \leq \int_{\partial \Omega} \log |p_n(z)| d\omega(z) = \int_{\partial \Omega \cap \D} \log |p_n(z)| d\omega(z) + \int_E \log |p_n(z)| d\omega(z).
\end{equation} 
The restrictions of the measures $\omega$ and $|dz|$ to $E$ are mutually absolutely continuous (this is a consequence of the rectifiability of $\partial \Omega$ and an old theorem of the Riesz brothers). Moreover, $\omega \leq |dz|$ on $E$ by the containment $\Omega \subset \D$ and the monotonicity of harmonic measure (\cite[Corollary 4.3.9]{ransford1995potential}). Using this, and an Egorov-type argument, one sees that the condition $\int_E |p_n(z)| |dz| \to 0$ implies that 
\[ \lim_{n \to \infty} \int_E \log |p_n(z)| d\omega(z) = - \infty.\] Thus, by \eqref{E:SubharmLogpn}, we will have $p_n(0) \to 0$ if we can control from above the quantity \[\int_{\partial \Omega \cap \D} \log |p_n(z)| d\omega(z).\] If $(p_n)_n$ obeys the growth condition imposed by the majorant $\lambda$, then a way to do so is to establish the integrability of $\lambda$ against the harmonic measure:
\begin{equation}
    \label{E:KhrushchevMajorantInt}
    \int_{\partial \Omega \cap \D} \log \lambda \bigg(\frac{1}{1-|z|} \bigg) d\omega(z) < +\infty.
\end{equation} Khrushchev does just that, and his domain of choice $\Omega$ is the unit disk from which he removed a sequence of curvlinear square boxes $S_\ell$, each having a complementary interval $\ell$ as one of the sides. See Figure~\ref{F:Figure3} (recall that in the classical \textit{Privalov cone domain} one removes instead triangular tents with bases $\ell$). A harmonic measure estimation in \cite[Lemma 3.1]{khrushchev1978problem} shows that on the relevant parts of the boundary of $S_\ell$, the harmonic measure $\omega$ is dominated by the arclength measure. Then the hypothesis \eqref{E:IntegralCarlesonCond} implies \eqref{E:KhrushchevMajorantInt}, and so $p_n(0) \to 0$. Since the harmonic measure $\omega(\Omega, \cdot, z)$ at points $z$ close to $0$ is comparable to $\omega(\Omega, \cdot, 0)$, we conclude that $p_n(z) \to 0$ in a neighbourhood of $0$. Khrushchev's theorem follows.

Observe, however, that the sides of the squares $S_\ell$ make a right angle with $\T$, and elementary formulas for conformal mappings onto domains with corners show that the harmonic measure near a corner should be significantly smaller than the arclength measure. In fact, we would expect \[\omega(I) \lesssim_{|\ell|} |I|^2\] whenever $I$ is a sufficiently short segment of $\partial S_\ell$ touching $\T$ (see Figure~\ref{F:Figure3}). Based on this reasoning, we may suspect the integrability condition \eqref{E:IntHypothesis} to be unnecessary, should the implied constant in the inequality above behave sufficiently well. We confirm this suspicion by instead considering a \textit{Joukowski-Privalov domain} $\D_E$, obtained from $\D$ by removing areas bounded by hyperbolic geodesics between endpoints of $\ell$ (see Figure~\ref{F:Figure2} below). Like Khrushchev's squares, the geodesics also meet $\T$ at a right angle, but in our new domain the harmonic measure may be easily and accurately estimated by explicit conformal mapping formulas involving the Joukowski function $z \mapsto z + 1/z$. Our main result in this direction is \thref{P:IntegrMajorantHarmMeas}, and may be of independent interest. In the proposition, we establish the corresponding integrability in \eqref{E:KhrushchevMajorantInt} for $\Omega = \D_E$ whenever $E$ is an $h$-Beurling-Carleson set and $h(t) = t \log \big(\lambda(1/t)\big)$, without any assumption regarding \eqref{E:IntHypothesis}. Together with Khrushchev's argument outlined above applied to $\Omega = \D_E$, \thref{T:MainTheorem1} is a direct corollary. 

\subsection{An application: sharpness of a local Jensen-type inequality} 
\label{S:ApplicationSubsec}
The earlier-mentioned weighted results from \cite{malman2023revisiting} have the following Fourier-analytic consequence (a proof is given in the preprint \cite{malman2023shift}).

\begin{thm*}
    Let $f$ be an integrable function on $\T$. Assume that we have the one-sided Fourier coefficient estimate 
    \begin{equation}
        \label{E:RapidDecay}
        |\widehat{f}(n)| = \mathcal{O}\big( \exp (-c \sqrt{n})\big), \quad n \geq 1
    \end{equation} for some $c > 0$. If $f$ does not vanish identically, then there exists an interval $I \subset \T$ for which
    \begin{equation}
        \label{E:fLogIntInterval}
        \int_I \log |f(z)| |dz| > -\infty.
    \end{equation} 
\end{thm*}

In fact, under the unilateral spectral decay \eqref{E:RapidDecay}, $\log |f|$ is locally integrable on the carrier $\{ z \in \T : |f(z)| > 0 \}$, a set which turns out to be essentially open (that is, differs from an open set only by a set of Lebesgue measure zero). The result is a type of \textit{local Jensen's inequality}, and it can be seen as a manifestation of a condensation phenomenon for functions obeying one-sided spectral estimates. Recall that the original \textit{global} Jensen's inequality says that a function $f$ with vanishing positive Fourier coefficients must satisfy $\log |f| \in L^1(|dz|)$, in fact $\log |\widehat{f}(0)| \leq \pi^{-1} \int_\T \log |f(z)| |dz|$. Improvements and variants of Jensen's inequality have a long history and appear notably in the works of Levinson \cite{levinson1940gap}, Beurling \cite{beurling1961quasianalyticity}, Volberg and Jöricke \cite{vol1987summability} and Borichev and Volberg \cite{borichev1990uniqueness}. The need for a local result came up in the study of the structure theory of de Branges' spaces $\hb$ (see \cite{malman2023shift}). 

The condition \eqref{E:RapidDecay} was a priori known to be close to optimal. Indeed, already Khrushchev's theorem stated above can be used to construct by duality a function $f$ living on a closed set $E$ containing no intervals, so that \eqref{E:fLogIntInterval} fails by default for any interval $I$, and yet $f$ satisfies 
\[ |\widehat{f}(n)| = \mathcal{O}\big( \exp (-c n^a)\big), \quad n \geq 1\] for any $a \in (0, 1/2)$. Our study of the potential for weakening of the hypotheses on the majorant in Khrushchev's Uniqueness Theorem was in part motivated by the desire to understand better the optimality of the condition \eqref{E:RapidDecay} in the context of detection of local logarithmic integrability of $|f|$. Our second main result shows that the condition is actually on point.

\begin{thm}
    \thlabel{T:LocalJensenSharpnessTheorem}
    Let $(c_n)_n$ be any sequence of positive real numbers which satisfies \[\lim_{n \to \infty} c_n = 0.\] There exists a non-zero Borel function $f$ on $\T$ for which the following two conditions hold:
    \begin{enumerate}[(i)]
        \item the set $E := \{ z \in \T : |f(z)| > 0 \}$ contains no intervals,
        \item we have the one-sided spectral estimate
    $|\widehat{f}(n)| = \mathcal{O}\big(\exp(-c_n \sqrt{n})\big)$, $n \geq 1$.
    \end{enumerate}
\end{thm}

We will see that \thref{T:LocalJensenSharpnessTheorem} is a consequence of \thref{T:MainTheorem1}. The carrying set $E$ of $f$ appearing in \thref{T:LocalJensenSharpnessTheorem} will be contained in an $h$-Beurling-Carleson set corresponding to $h(x)$ with sufficiently slow decay as $x \to 0+$.

\subsection{Outline of the rest of the note}
In Section~\ref{S:Section2} we estimate the harmonic measure on parts hyperbolic geodesics in the upper half-plane. Section~\ref{S:Section3} deals with the proof of integrability of the majorant against the harmonic measure of a Joukowski-Privalov type domain. In the last Section~\ref{S:Section4} we show how \thref{T:LocalJensenSharpnessTheorem} follows from \thref{T:MainTheorem1}.

\section{Harmonic measure on hyperbolic geodesics}
\label{S:Section2}

\subsection{A Joukowski domain}

Let $\H = \{ z \in \mathbb{C} : \Im z > 0 \}$ denote the upper half-plane. The domain 
\begin{equation}
    \label{E:OmegaLDef}
    \Omega_L := \H \setminus \{ z \in \mathbb{C} : |z| \leq L \}
\end{equation} has as its boundary the two line segments $[-\infty, -L]$, $[L, \infty]$ and the hyperbolic geodesic (a half-circle) between $-L$ and $L$: \[A_L := \{ z \in \mathbb{C} : |z| = L, \Im z > 0 \}.\] The geodesic meets the real axis in two orthogonal angles. We shall see that the Joukowski function $z \mapsto z + 1/z$ can be used to express the conformal mappings from $\Omega_L$ onto the upper half-plane $\H$. See \cite[Chapter 6]{kythe2019handbook} for similar maps involving the Joukowski function.

\begin{figure}
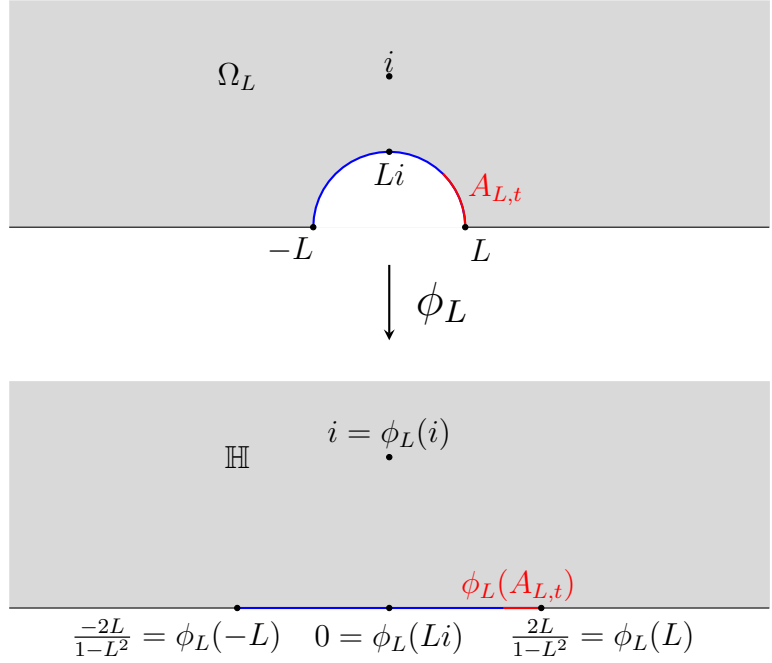

    \centering
    \includestandalone{joukowski}
    \hfill     
    \includestandalone{HalfplaneMarked}
    \caption{The action of the mapping $\phi_L$ between $\Omega_L$ and $\mathbb{H}$.}
    \label{F:Figure1}
\end{figure}

In our presentation, we shall assume that $L$ is small, so that we may suppose $i \in \Omega_L$. An elementary computation shows that the mapping 
\[ \phi_L(z) := \frac{L}{1-L^2}\Bigg(\frac{L}{z} + \frac{z}{L} \Bigg), \quad z \in \Omega_L\] is a conformal bijection from $\Omega_L$ onto $\H$ which satisfies \[ \phi_L(i) = i.\] To see that the map is conformal between the domains, consider for a moment $\phi_L$ as a mapping with domain of definition $\mathbb{C} \setminus \{ 0 \}$ and note that its imaginary part is 
\[ \Im \phi_L(re^{it}) = \frac{L}{1-L^2}\Bigg( -\frac{L}{r} + \frac{r}{L}\Bigg) \sin(t), \quad re^{it} \in \mathbb{C} \setminus \{ 0 \}.\]  We check easily that this quantity is positive precisely when $z = re^{it}$ lies in $\Omega_L$ or in the lower half of the disk of radius $L$ centered at the origin. That is 
\begin{equation}
    \label{E:EquivPhiLOmegaL}
    \phi_L(z) \in \H \quad \Longleftrightarrow \quad z \in \Omega_L \cup \{ z \in \mathbb{C} : |z| < L, \Im z < 0 \}.
\end{equation} If $w \in \H$, then solving for $z$ in the equation $\phi_L(z) = w$ leads to a quadratic expression with two solutions $z_1, z_2$ which, by the equivalence in \eqref{E:EquivPhiLOmegaL}, must both lie in the union of the two domains appearing in \eqref{E:EquivPhiLOmegaL}. The relation $\phi_L(z) = \phi_L(\frac{L^2}{z})$ shows that the involutive mapping $z \mapsto \frac{L^2}{z}$ interchanges the two solutions. The mapping also interchanges the two domains $\Omega_L$ and $\{ z \in \mathbb{C} : |z| < L, \Im z < 0 \}$. This shows that exactly one of the solutions $z_1, z_2$ lies in $\Omega_L$. We conclude that $\phi_L: \Omega_L \to \H$ is bijective. See Figure~\ref{F:Figure1} for an illustration of the action of the mapping $\phi_L$.

\subsection{Harmonic measure of subarcs}

Consider the arc
\[ A_{L,t} = \{ Le^{i\theta} : 0 \leq \theta \leq t \}, \quad t < \pi/2.\] 
This is a subarc of $A_L$ of length $|A_{L,t}| = Lt$ lying in the semicircular part of the boundary of $\Omega_L$, with right end-point at $L \in \R \cap \partial \Omega_L$. See the red markings in Figure~\ref{F:Figure1}. We are interested in the harmonic measure of this arc with respect to the point $i \in \Omega_L$. Since the arc $A_{L,t}$ lies near a corner of $\Omega_L$, the harmonic measure $\omega(\Omega_L, A_{L,t}, i)$ should be significantly smaller than the length of the arc. We will prove the following estimate.

\begin{lem}
    \thlabel{L:ArcHarmMeasEst}
    In the notation as above, we have
    \[ \omega(\Omega_L, A_{L,t}, i) \leq (1-\cos(t))\frac{2L}{\pi(1-L^2)}.\]
\end{lem}

\begin{proof}
Since the mapping $\phi_L$ introduced above is conformal and satisfies $\phi_L(i) = i$, we know by invariance of harmonic measure that 
\[ \omega(\Omega_L, A_{L,t}, i) = \omega( \mathbb{H}, \phi_L(A_{L,t}), i).\] Here
\[ \phi_L(A_{L,t}) = \{x \in \partial \mathbb{H} = \R : x = \phi_L(z), z \in A_{L,t} \}\] is an interval in $\R$ with right endpoint at $\phi_L(L) = \frac{2L}{1 - L^2}$. The left endpoint is given by 
\[ \phi_L(Le^{it}) = \frac{L}{1-L^2}\bigg( \frac{Le^{it}}{L} + \frac{L}{Le^{it}}\bigg) = \frac{2L \cos (t)}{1-L^2}.\] The well-known upper half-plane harmonic measure formula
\[ \omega( \mathbb{H}, B, i) = \frac{1}{\pi}\int_B \frac{1}{t^2+1} \, dt\] which holds for Borel subsets $B \subset \R$, readily implies the desired estimate:
\begin{align*}
    \omega(\Omega_L, A_{L,t}, i) &= \frac{1}{\pi}\int^{\frac{2L}{1-L^2}}_{\frac{2L \cos (t)}{1-L^2}} \frac{1}{t^2 + 1} \, dt \\
    &\leq \frac{1}{\pi}\int^{\frac{2L}{1-L^2}}_{\frac{2L \cos (t)}{1-L^2}} 1 \, dt \\ 
    &= (1-\cos(t))\frac{2L}{\pi(1-L^2)}.    
\end{align*} 
\end{proof}

\section{Integrability properties of harmonic measure in Joukowski-Privalov domains}

\label{S:Section3}

\subsection{A Joukowski-Privalov domain} \label{S:JoukowskiPrivalovDefSubsec} 

\begin{figure}
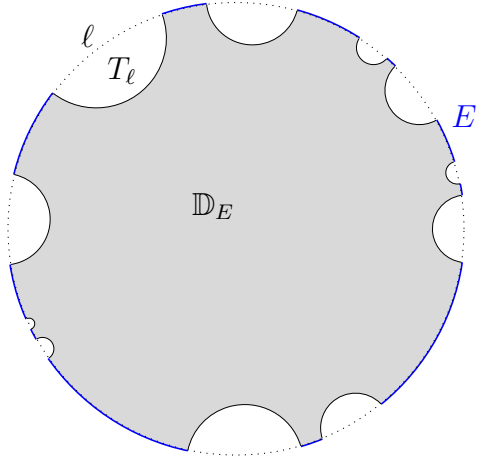

    \centering
    \includestandalone{Joukowski-Privalov}
    \caption{An example of a domain $\D_E$, with the set $E$ marked in blue.}
    \label{F:Figure2}
\end{figure}

Let $E$ be a closed subset of $\T$ and $\C(E)$ be the family of maximal open intervals complementary to $E$ on $\T$. By adding a finite number of points to $E$, we may assume that all intervals in $\C(E)$ have length bounded above by any small a priori chosen constant (such a transformation of $E$ does not change its membership in any of the $h$-Beurling-Carleson classes, and it does not change its measure, so all our results can be translated back and forth between any two sets differing by a finite number of points). Then to every $\ell \in \C(E)$ corresponds a unique hyperbolic geodesic $\gamma_\ell \subset \D$ which connects the two endpoints of $\ell$. The curve $\gamma_\ell$ is a circular arc which meets $\T$ in two orthogonal angles. The interval $\ell$ and the curve $\gamma_\ell$ enclose a small region $T_\ell \subset \D$. We define the domain $\D_E$ as the difference between the unit disk $\D$ and the union of closures of $T_\ell$, namely \[ \D_E := \D - \big( \bigcup_{\ell \in \C(E)} \conj{T_\ell} \big). \] 
See Figure~\ref{F:Figure2} for a picture of an example domain $\D_E$. 

For each $e^{it} \in \T$, there exist a largest $r = r(t) \in (0, 1]$ for which we have $re^{it} \in \partial \D_E$. The mapping $e^{it} \to re^{it}$ is easily seen to be continuous, so $\partial \D_E$ is a Jordan curve. The lengths of $\gamma_\ell$ and $\ell$ are comparable, and it follows easily that $\partial \D_E$ is a rectifiable curve. 

\subsection{Majorant integrated against the harmonic measure}

Let $h$ be an increasing function which defines a majorant as in \eqref{E:LessThanExpGrowth}. We assume that $h$ satisfies the regularity hypotheses stated in \eqref{E:RegCond} in Section~\ref{S:MainResultSec}.

\begin{prop}
    \thlabel{P:IntegrMajorantHarmMeas}
    Let \[\omega = \omega(\D_E, \cdot, 0)\] denote the harmonic measure at $0$ for the domain $\D_E$ introduced above. Assume that $E$ is a subset of $\T$ which is closed, and satisfies the $h$-Beurling-Carleson condition in \eqref{E:CarlesonCondDef}. Then we have
    \[ \int_{\partial \D_E \cap \D} \frac{h(1-|z|^2)}{1-|z|^2} d\omega(z) < \infty.\]
\end{prop}

Note that \eqref{E:RegCond} implies 
\[ \frac{h(1-|z|^2)}{1-|z|^2} \leq \frac{h(1-|z|)}{1-|z|} \leq 2\frac{h(1-|z|^2)}{1-|z|^2}. \]

To prove \thref{P:IntegrMajorantHarmMeas}, we will use the well-known subordination principle for harmonic measure.

\begin{lem}
    \thlabel{L:SubOrdHarmMeasure}
    Let $\ell \in \C(E)$ be one of the intervals complementary to $E$. Denote by $\ell^c = \T \setminus \ell$ its closed complement in $\T$, and by $\D_{\ell^c} = \D \setminus T_\ell$ the corresponding Joukowski-Privalov domain which has only one removed region. If $B$ is a Borel subset of $\gamma_\ell \subset \partial \D_E \cap \partial \D_{\ell^c}$, then 
    \[ \omega(B) = \omega(\D_E, B, 0) \leq \omega(\D_{\ell^c}, B, 0).\]
\end{lem}

The lemma is a consequence of the containment $\D_E \subset \D_{\ell^c}$ and the maximum principle for harmonic functions. We refer to \cite[Corollary 4.3.9]{ransford1995potential} for a proof.

We will also use a very simple distortion estimate. 

\begin{lem}
    \thlabel{L:DistEst}
    Let $\psi(z) = i \frac{1-z}{1+z}$, $z \in \D$. For every $z_1,z_2 \in \D$, $\Re z_1 \geq 0$, $\Re z_2 \geq 0$, we have the inequalities
    \[ \frac{|z_1 - z_2|}{2} \leq |\psi(z_1) - \psi(z_2)| \leq 2|z_1 - z_2|\]
\end{lem}

The result follows easily by a direct computation and/or elementary estimation of the derivative of $\phi$ and its inverse. We leave out a proof.


\begin{proof}[Proof of \thref{P:IntegrMajorantHarmMeas}]
    Since 
    \[ \int_{\partial \D_E \cap \D} \frac{h(1-|z|^2)}{1-|z|^2} d\omega(z) = \sum_{\ell \in \C(E)} \int_{\gamma_\ell} \frac{h(1-|z|^2)}{1-|z|^2} d\omega(z),\] by our assumption of finiteness of $\sum_{\ell \in \C(E)} h(|\ell|)$ it will suffice to show that
    \begin{equation}
        \label{E:SuffIntEstimate}
        \int_{\gamma_\ell} \frac{h(1-|z|^2)}{1-|z|^2} d\omega(z) \leq C h(|\ell|) 
    \end{equation} where the constant $C > 0$ does not depend on the complementary interval $\ell$. Fix one of the intervals $\ell$ for the remainder of the proof. Using \thref{L:SubOrdHarmMeasure} we conclude that 
    \begin{equation}
        \label{E:FirstEstimate}
        \int_{\gamma_\ell} \frac{h(1-|z|^2)}{1-|z|^2} d\omega(z) \leq \int_{\gamma_\ell} \frac{h(1-|z|^2)}{1-|z|^2} d\omega_\ell(z) 
    \end{equation} where $\omega_\ell = \omega(\D_{\ell^c}, \cdot, 0)$ denotes the harmonic measure in $\D_{\ell^c}$. By a rotation of the domain $\D_{\ell^c}$ and the radial symmetry of the integrand above, may suppose that the interval $\ell$ has $z = 1$ as its midpoint. The conformal mapping \[ \psi(z) = i \frac{1-z}{1+z}, \quad z \in \D,\] which we considered in \thref{L:DistEst}, maps the domain $\D_{\ell^c}$ bijectively onto a subdomain of the upper half-plane $\H$, and it satisfies $\psi(0) = i$, $\psi(1) = 0$. Since $\psi$ is a linear fractional transformation, it transforms the circular arc $\gamma_\ell$ into a circular arc $\psi(\gamma_\ell)$ which is symmetric with respect to the imaginary axis, and meets the real axis in two orthogonal angles. We conclude that $\psi(\D_{\ell^c})$ is one of the earlier introduced domains $\Omega_L$ (recall the definition in \eqref{E:OmegaLDef}). The interval $[-L, L]$ is the image of $\ell$ under $\psi$. Since we made arrangements for $\ell$ to be short, the arc $\gamma_\ell$, and the domain it encloses together with $\ell$, are located well within the right half-plane. Therefore, by the distortion estimates in \thref{L:DistEst}, the lengths of the intervals $[-L,L]$ and $\ell$ are comparable. More precisely, we have 
    \begin{equation}
        \label{E:ellLcompare}
        |\ell|/2 \leq 2L \leq 2|\ell|.
    \end{equation} Changing variables conformally, we obtain
    \begin{equation}
        \label{E:VariableChange}
        \int_{\gamma_\ell} \frac{h(1-|z|^2)}{1-|z|^2)} d\omega_\ell(z) = \int_{A_L} \frac{h(1-|\psi^{-1}(w)|^2)}{1-|\psi^{-1}(w)|^2} d\omega_L(w) 
    \end{equation} where now \[\omega_L := \omega(\Omega_L, \cdot, i)\] is the harmonic measure in $\Omega_L$, and as before,\[ A_L = \{ z \in \H : |z| = L, \Im z > 0 \}\] is the circular part of the boundary of $\Omega_L$. Decompose the part of $A_L$ lying in the right half-plane as 
    \[ A_L \cap \{ z \in \mathbb{C} : \Re z \geq 0\} = \bigcup_{n=0}^\infty B_{L,n}\] where 
    \[ B_{L,n} = \big\{ z = Le^{it} \in A_L : t_{n+1} < t \leq t_n \big\}\] and
    \[ t_n := \frac{\pi}{2}2^{-n}, \quad n \geq 0.\]
    Since $B_{L,n} \subset A_{L,t_n}$, we obtain by \thref{L:ArcHarmMeasEst} the estimation
    \begin{align*}
        \omega_L(B_{L,n}) &\leq (1-\cos(t_n)) \frac{2L}{\pi(1-L^2)} \leq t_n^2 L
    \end{align*} where we have used the inequality \[ 1-\cos(t) \leq t^2/2, \quad t > 0\] and implicitly assumed that $|\ell|$, and therefore $L$, are sufficiently small to have $\frac{1}{\pi(1-L^2)} \leq 1$ (recall the discussion at the beginning of Section~\ref{S:JoukowskiPrivalovDefSubsec}).

    Let us now compute 
    \[ \sup_{w \in B_{L,n}} \frac{h(1-|\psi^{-1}(w)|^2)}{1-|\psi^{-1}(w)|^2}.\]
    Since 
    \[ \psi^{-1}(w) = \frac{i-w}{i+w}\]
    we obtain, for $w \in B_{L,n}$, that
      \[  1- |\psi^{-1}(w)|^2 = \frac{4 \Im w}{|i+w|^2} \geq \Im w.\] The inequality above is a consequence of
      \[ |i+w| \leq 1+|w| = 1+L < 2\] (recall that we arranged for $i \in \Omega_L$, so $L < 1$). Using that $\sin(t) \geq t/2$ for $t \in [0, \pi/2]$, we obtain \[\inf_{w \in B_{L,n}} \Im w = L\sin(t_{n+1}) \geq L \frac{t_{t+1}}{2} = L \frac{t_n}{4}. \] Thus
      \[ \inf_{w \in B_{L,n}} 1- |\psi^{-1}(w)|^2 \geq L \frac{t_n}{4}, \] and since $h(x)/x$ is assumed decreasing by \eqref{E:RegCond}, we obtain finally that
      \[ \sup_{w \in B_{L,n}} \frac{h(1-|\psi^{-1}(w)|^2)}{1-|\psi^{-1}(w)|^2} \leq \frac{h(Lt_n/4)}{Lt_n/4}.\] By combining this estimate with our earlier estimate for $\omega_L(B_{L,n})$, it follows that
      \[ \int_{B_{L,n}} \frac{h(1-|\psi^{-1}(w)|^2)}{1-|\psi^{-1}(w)|^2} d\omega_L(w) \leq 4 h(Lt_n/4) t_n \leq 2 \pi h(|\ell|) 2^{-n}.\] In the last inequality we used that $h$ is increasing and that $Lt_n/4 \leq L$ for every $n \geq 0$ (recall \eqref{E:ellLcompare}). Summing over $n$, and employing a similar argument for the part of $A_L$ contained in the left half-plane, we obtain
      \[ \int_{A_L} \frac{h(1-|\psi^{-1}(w)|^2)}{1-|\psi^{-1}(w)|^2} d\omega_L(w) \leq 8 \pi h(|\ell|).\]
      This last estimate combined with \eqref{E:SuffIntEstimate}, \eqref{E:FirstEstimate} and \eqref{E:VariableChange} proves the desired result. 
\end{proof}

As explained in Section~\ref{S:TechniqueSection}, Khrushchev's technique combined with \thref{P:IntegrMajorantHarmMeas} immediately implies \thref{T:MainTheorem1}.

\section{Sharpness of the local Jensen-type inequality}
\label{S:Section4}

Here, we prove \thref{T:LocalJensenSharpnessTheorem}.

\subsection{Lemmas on estimation of moments}

To establish the theorem, we will need to go through a rather tedious but fairly standard sequence of lemmas related to estimation of moments of weights related to the majorants in \eqref{E:LessThanExpGrowth}.

Note first that \thref{T:LocalJensenSharpnessTheorem} holds for $(c_n)_n$ if it holds for some larger sequence $(\widetilde{c}_n)_n$, i.e, one satisfying $c_n \leq \widetilde{c}_n$ for all $n$. Thus we may replace $(c_n)_n$ by any larger sequence which satisfies some additional useful regularity properties.

\begin{lem} \thlabel{L:SequenceReplacementLemma}
    Let $(c_n)_n$ be a sequence of positive reals which satisfies $\lim_{n \to \infty} c_n = 0$. There exists a sequence of positive reals $(\widetilde{c}_n)_n$ that satisfies $\lim_{n \to \infty} \widetilde{c}_n = 0$ and the following three additional conditions:
    \begin{enumerate}[(i)]
        \item $c_n \leq \widetilde{c}_n$ for all $n$,
        \item $(\widetilde{c}_n)_n$ is decreasing,
        \item $\widetilde{c}_n \leq \sqrt{\frac{n+1}{n}} \widetilde{c}_{n+1}$ for all $n$.
    \end{enumerate}
\end{lem}

\begin{proof}
We define the desired sequence recursively by setting $\widetilde{c}_1 = c_1$ and then 
\[ \widetilde{c}_{n+1} = \max\Bigg(  \max_{m \geq n+1} c_m , \sqrt{\frac{n}{n+1}} \widetilde{c}_n \Bigg).\]
A minute of reflection shows that the conditions $(i)$, $(ii)$ and $(iii)$ hold for $(\widetilde{c}_n)_n$. To see that this sequence decays to $0$, we use $(ii)$ to conclude that if $\widetilde{c}_{n+1} = \max_{m \geq n+1} c_m$ for infinitely many $n$, then the claim holds, since $\lim_{n \to \infty} \max_{m \geq n+1} c_m = 0$ by hypothesis. In the contrary case, we have $\widetilde{c}_{n+1} = \sqrt{\frac{n}{n+1}}\widetilde{c}_n$ for all $n \geq n_0$, say. In that case, we obtain recursively that
\[ \widetilde{c}_{n_0+m} \leq \sqrt{\frac{n_0}{n_0+m}} \widetilde{c}_{n_0}\] and so $\widetilde{c}_{n_0+m} \to 0$ as $m \to \infty$.    
\end{proof}

An elementary calculus argument, which we skip, establishes the following formula.

\begin{lem} \thlabel{L:LegendreTransformOneOverX}
Let $c > 0$. For large positive integers $n$, we have
\[ \inf_{x \in (0, 1)} nx + \frac{c^2}{x} = 2c \sqrt{n},\] and the infimum is attained at the point $x_n = \frac{c}{\sqrt{n}} \in (0,1)$.    
\end{lem}

By \thref{L:SequenceReplacementLemma}, we may assume that $(c_n)_n$ satisfies the three conditions appearing in the statement of the lemma. Clearly, the sequence $\big( \frac{c_n}{\sqrt{n}}\big)_{n \geq 1}$ is decreasing to $0$, and we may assume that $c_1 < 1$. Define a function $h(x)$ by
\begin{equation}
    \label{E:hDefForCseq}
    h(x) = c_n^2, \quad x \in \Big( \frac{c_{n+1}}{\sqrt{n+1}}, \frac{c_{n}}{\sqrt{n}}\Big].
\end{equation}
For definiteness, we may set $h(x) = c_1^2$ for $x > c_1$. 

\begin{lem} \thlabel{L:MainEstimate}
Let $h$ be defined as above. For large positive integers $n$, we have 
\[ \inf_{x \in (0, 1)} nx + \frac{h(x)}{x} \geq c_n\sqrt{n}.\]
\end{lem}

\begin{proof}
    Consider first the interval \[ I_n := [ c_n/\sqrt{n}, 1]. \] Since $(c_n)_n$ is decreasing, by our definition of $h$ we have the inequality
    \[ c_n^2 \leq h(x), \quad x \in I_n.\]
    Using \thref{L:LegendreTransformOneOverX}, we obtain
    \begin{align*}
        \inf_{x \in I_n} nx + \frac{h(x)}{x} \geq \inf_{x \in I_n} nx + \frac{c_n^2}{x} \geq c_n\sqrt{n}.
    \end{align*} Now consider the interval \[J_m := \Big( \frac{c_{m+1}}{\sqrt{m+1}}, \frac{c_{m}}{\sqrt{m}}\Big], \quad m \geq n.\] By \eqref{E:hDefForCseq}, $h \equiv c_m^2$ on $J_m$. An elementary calculus argument reveals that $nx + \frac{h(x)}{x}$ attains its infimum over $J_m$ at the right endpoint of the interval. Indeed, in the interior of $J_m$, the derivative of the expression equals
    \[ n - \frac{c_m^2}{x^2}\] and this quantity is negative if \[ x < \frac{c_m}{\sqrt{n}},\] so certainly it is negative for $x$ in the interior of $J_m$, since $m \geq n$. We conclude that $x \mapsto nx + \frac{h(x)}{x}$ is decreasing in the interval $J_m$ and thus
    \[ \inf_{x \in J_m} nx + \frac{h(x)}{x} = n \frac{c_m}{\sqrt{m}} + c_m \sqrt{m}.\] We claim that
    \[ n \frac{c_m}{\sqrt{m}} + c_m \sqrt{m} \geq c_n \sqrt{n}.\] If this inequality holds, then our proof is complete. Dividing through by $\sqrt{n}$, we see that the inequality is equivalent to
    \[ \bigg(\sqrt{\frac{n}{m}} + \sqrt{\frac{m}{n}} \bigg) c_m \geq c_n.\] The inequality will hold if the stronger inequality
    \[ \sqrt{\frac{m}{n}} c_m \geq c_n\] is satisfied. Since $m \geq n$, this last inequality just the iterated version of part $(iii)$ of \thref{L:SequenceReplacementLemma}. 
\end{proof}

The function $h$ defined above is not continuous. Replacing $h$ by its least concave majorant, we conclude that given any sequence of positive reals $(c_n)_n$ which satisfies $\lim_{n \to \infty} c_n = 0$, there exists a continuous function $h$ satisfying $h(0) = 0$ and the regularity hypotheses in \eqref{E:RegCond} of Section~\ref{S:MainResultSec} for which the inequality 
\begin{equation}
    \label{E:hLegendreDominateMoments} \inf_{x \in (0,1)} nx + \frac{h(x)}{x} \geq c_n \sqrt{n}
\end{equation} holds for large positive integers $n$.

\begin{lem} \thlabel{L:MomentEstimate}
Assume that $h$ satisfies \eqref{E:hLegendreDominateMoments}. Then for large positive integers $n$ we have the inequality
\[ \int_0^1 x^n \exp\bigg(-\frac{h(1-x)}{1-x}\bigg) \, dx \leq \exp(-c_n\sqrt{n}).\]
\end{lem}

\begin{proof}
    Since the integration interval is of length $1$, it suffices to establish that the integrand is pointwise bounded by $\exp( -c_n \sqrt{n})$. To do so, we use the inequality $\log x \leq x-1$ which holds for all $x > 0$, and the change of variable $1-x = s$:
    \begin{align*}
        x^n \exp\bigg(-\frac{h(1-x)}{1-x}\bigg) &= \exp\bigg(n \log x  -\frac{h(1-x)}{1-x}\bigg) \\ 
        &\leq \exp\bigg(n(x-1) -\frac{h(1-x)}{1-x}\bigg) \\
        &=  \exp\bigg( - ns - \frac{h(s)}{s}\bigg) \\
        &\leq \exp( -c_n \sqrt{n}).
    \end{align*}
    In the last step we used \eqref{E:hLegendreDominateMoments}.
\end{proof}

\subsection{Proof of \thref{T:LocalJensenSharpnessTheorem}}
Let $(c_n)_n$ be any sequence satisfying $\lim_{n \to \infty} c_n = 0$ and $h$ be a corresponding function satisfying \eqref{E:RegCond} and the inequalities in \eqref{E:hLegendreDominateMoments}. Let \[ G(x) = \exp\bigg(-\frac{h(1-x)}{1-x} \bigg), \quad x \in [0,1])\] and $\Po^2_G$ be the weighted Bergman space consisting of analytic functions in $\D$ that satisfy
\[ \| f\|^2_G := \int_\D |f(z)|^2 G(|z|) dA(|z|) < \infty.\] The space $\Po^2_G$ is the closure of analytic polynomials in the corresponding weighted Lebesgue space $L^2(G\, dA)$. A short computation involving a subharmonicity argument and \eqref{E:RegCond} shows that for $z \in \D$ we have 
\begin{align*}
    |f(z)|^2 &\leq  \frac{4}{\pi (1-|z|)^2} \int_{ w: |z-w| < (1-|z|)/2} |f(w)|^2 dA(w) \\
    &\leq \frac{C_1 \|f\|^2_G}{(1-|z|)^2} \exp \bigg( C_2 \frac{h(1-|z|)}{1-|z|} \bigg) 
\end{align*}
i.e, any bounded sequence in $\Po^2_G$ is controlled pointwise in $\D$ by a majorant $\lambda_h$ from the class \eqref{E:LessThanExpGrowth} (after a replacement $C_2 h \to h$) if $h$ decays slow enough to subsume the term $(1-|z|)^{-2} = \exp \big(- 2\log(1-|z|) \big)$ into the exponent by increasing $C_2$ slightly. Should $h$ decay too fast we may always replace $h$ by a larger function which dominates $x \log(1/x)$ and for which the above subsumption is valid (for instance, we may replace $h(x)$ by $h(x) + x\log(1/x)$). The inequality \eqref{E:hLegendreDominateMoments} still holds for this larger $h$.

Let $E$ be an $h$-Beurling-Carleson set of positive measure which contains no intervals, and consider the diagonal 
\[ \mathcal{D} := \{ p \in \Po : (p,p) \} \subset \Po^2_G \oplus L^1(E),\] the second factor being the usual Lebesgue space on $E \subset \T$, defined using the arclength measure $|dz|$ restricted to $E$. The diagonal $\mathcal{D}$ is not norm-dense in the product space. Indeed, by \thref{T:MainTheorem1} and the pointwise estimate above, no tuple of the form $(F, 0)$, $F \neq 0$ is a norm-limit of a sequence in $\mathcal{D}$. A non-zero functional therefore annihilates $\mathcal{D}$. This functional is represented in the usual way by a tuple of functions $g \in \Po^2_G$, $f \in L^\infty(E)$. Since $(z^n, z^n) \in \mathcal{D}$ and $f$ lives only on $E$, we have
\[ \widehat{f}(n) = \int_E f(z) \conj{z}^n \, |dz| = -\int_\D g(z)\conj{z^n} G(|z|) dA(z), \quad n \geq 0.\] The function $f$ is non-zero: were it zero, then the above relation shows that Taylor coefficients around the origin of the analytic function $g$ vanish, and then so does $g$, which is contrary to hypothesis of the functional being non-zero. Cauchy-Schwarz inequality applied to $|g(z)|\sqrt{|z|^nG(|z|)}$ and $\sqrt{|z|^nG(|z|)}$, and the moment estimate in \thref{L:MomentEstimate} show that
\[ |\widehat{f}(n)| \leq \|g\|_G \cdot \exp\big(-(c_n/2) \sqrt{n} \big).\] Since $(c_n)_n$ is an arbitrary positive sequence satisfying $\lim_{n \to \infty} c_n = 0$, our \thref{T:LocalJensenSharpnessTheorem} follows.

\bibliographystyle{alpha}
\bibliography{mybib}

\end{document}

%% file: khrushchevsquare.tex
\begin{tikzpicture}[scale=3]    
    \fill[gray, opacity=0.3] (0,1) arc(90:60:1) -- ({0.738*cos(60)}, {0.738*sin(60)}) arc(60:45:0.738) -- ({cos(45)}, {sin(45)}) arc(45:30:1) --  ({0.869*cos(30)}, {0.869*sin(30)}) arc(30:22.5:0.869) -- ({cos(22.5)}, {sin(22.5)}) arc(22.5:0:1) -- (0,0);

    \draw[very thin] (0,1) arc(90:60:1) -- ({0.738*cos(60)}, {0.738*sin(60)}) arc(60:45:0.738) -- ({cos(45)}, {sin(45)}) arc(45:30:1) --  ({0.869*cos(30)}, {0.869*sin(30)}) arc(30:22.5:0.869) -- ({cos(22.5)}, {sin(22.5)}) arc(22.5:0:1);

    \draw[dotted, very thin] (0,1) arc(90:0:1);

    \node at (0.2, 0.2) {$\Omega$};
    \node at (0.50, 0.70) {\tiny $S_\ell$};
    \node at (0.6, 0.85) {\tiny $\ell$};
    \node[color=blue] at (0.70, 0.62) {\tiny $I$};

    \draw[thick, color=blue] ({cos(45)}, {sin(45)}) -- ({0.85*cos(45)}, {0.85*sin(45)});

\end{tikzpicture}

%% file: joukowski.tex
\begin{tikzpicture}[scale=1]    
    \filldraw[gray, opacity=0.3] (-5, 0) -- (5, 0) -- (5, 3) -- (-5, 3);
    \filldraw[white] (-1,0) -- (1,0) arc(0:180:1);
    \draw [blue,thick, domain=0:180] plot ({cos(\x)}, {sin(\x)});

    \draw [red,thick, domain=0:45] plot ({cos(\x)}, {sin(\x)});
    \node[red] at (1.4, 0.5) {$A_{L,t}$};

    \draw [black, thin] (-5,0) -- (-1,0);
    \draw [black, thin] (1,0) -- (5,0);
    \draw [white, thin] (1,-0.2) -- (5,-0.2);
   
    \node at (-2, 2) {$\Omega_L$};

    \filldraw (0, 2) circle (1pt); 
    \node at (0, 2.2) {$i$};

    \filldraw (0, 1) circle (1pt); 
    \node at (0, 0.7) {$Li$};

    \filldraw (1, 0) circle (1pt); 
    \node at (1.2, -0.3) {$L$};
     	    
    \filldraw (-1, 0) circle (1pt); 
    \node at (-1.3, -0.3) {$-L$};

    \draw [-stealth, thick] (0, -0.5) -- (0, -1.5);
    \node at (0.7, -1) {\scalebox{1.5}{$\phi_L$}};

    \draw [white, thin] (1,-2) -- (5,-2);

\end{tikzpicture}

%% file: HalfplaneMarked.tex
\begin{tikzpicture}[scale=1]    
    \filldraw[gray, opacity=0.3] (-5, 0) -- (5, 0) -- (5, 3) -- (-5, 3);
    \draw [black, thin] (-5,0) -- (5,0);
    \draw [blue, thick] (-2,0) -- (2,0);
    \draw [red, thick] (1.5,0) -- (2,0);
    \draw [white, thin] (1,-0.2) -- (5,-0.2);
    \node[red] at (1.7, 0.3) {$\phi_L(A_{L,t})$};bottom white strip
   
    \node at (-2, 2) {$\mathbb{H}$};

    \filldraw (0, 2) circle (1pt); 
    \node at (0, 2.3) {$i = \phi_L(i)$};

    \filldraw (0, 0) circle (1pt); 
    \node at (0, -0.4) {$0 = \phi_L(Li)$};

    \filldraw (2, 0) circle (1pt); 
    \node at (2.8, -0.4) {$\frac{2L}{1-L^2} = \phi_L(L)$};
     	    
    \filldraw (-2, 0) circle (1pt); 
    \node at (-2.8, -0.4) {$\frac{-2L}{1-L^2} = \phi_L(-L)$};

\end{tikzpicture}

%% file: Joukowski-Privalov.tex
\begin{tikzpicture}[scale=3]    
  \tkzDefPoint(0,0){O}
  \tkzDefPoint(1,0){A}
  \tkzDrawCircle[color=blue, very thick](O,A)
  \tkzFillCircle[color=gray!30](O,A)

  \node at (-0.1, 0.1) {$\mathbb{D}_E$};
  \node[blue] at (1, 0.5) {$E$};

  \tkzDefPoint(cos(0.2), sin(0.2)){z1}
  \tkzDefPoint(cos(0.3), sin(0.3)){w1} 
  \tkzDefCircle[orthogonal through=z1 and w1](O,A) \tkzGetPoint{B1}
  \tkzInterCC(O,A)(B1,z1) \tkzGetPoints{C1}{D1}   
  \tkzFillCircle[color=white](B1,z1)  
  \tkzDrawArc(B1,C1)(D1)

  \tkzDefPoint(.9,.1){z2}
  \tkzDefPoint(.9,-.1){w2} 
  \tkzDefCircle[orthogonal through=z2 and w2](O,A) \tkzGetPoint{B2}
  \tkzInterCC(O,A)(B2,w2) \tkzGetPoints{C2}{D2}     
  \tkzFillCircle[color=white](B2,z2)
  \tkzDrawArc(B2,C2)(D2)

  \tkzDefPoint(cos(0.5), sin(0.5)){z3}
  \tkzDefPoint(cos(0.8), sin(0.8)){w3} 
  \tkzDefCircle[orthogonal through=z3 and w3](O,A) \tkzGetPoint{B3}
  \tkzInterCC(O,A)(B3,w3) \tkzGetPoints{C3}{D3}   
  \tkzFillCircle[color=white](B3,w3)  
  \tkzDrawArc(B3,C3)(D3)

  \tkzDefPoint(cos(0.85), sin(0.85)){z4}
  \tkzDefPoint(cos(1), sin(1)){w4} 
  \tkzDefCircle[orthogonal through=z4 and w4](O,A) \tkzGetPoint{B4}
  \tkzInterCC(O,A)(B4,w4) \tkzGetPoints{C4}{D4}   
  \tkzFillCircle[color=white](B4,w4)  
  \tkzDrawArc(B4,C4)(D4)

  \tkzDefPoint(cos(1.3), sin(1.3)){z5}
  \tkzDefPoint(cos(1.7), sin(1.7)){w5} 
  \tkzDefCircle[orthogonal through=z5 and w5](O,A) \tkzGetPoint{B5}
  \tkzInterCC(O,A)(B5,w5) \tkzGetPoints{C5}{D5}   
  \tkzFillCircle[color=white](B5,w5)  
  \tkzDrawArc(B5,C5)(D5)      

  \tkzDefPoint(cos(1.9), sin(1.9)){z6}
  \tkzDefPoint(cos(2.5), sin(2.5)){w6} 
  \tkzDefCircle[orthogonal through=z6 and w6](O,A) \tkzGetPoint{B6}
  \tkzInterCC(O,A)(B6,w6) \tkzGetPoints{C6}{D6}   
  \tkzFillCircle[color=white](B6,w6)  
  \tkzDrawArc(B6,C6)(D6)    

  \tkzDefPoint(cos(2.9), sin(2.9)){z7}
  \tkzDefPoint(cos(3.3), sin(3.3)){w7} 
  \tkzDefCircle[orthogonal through=z7 and w7](O,A) \tkzGetPoint{B7}
  \tkzInterCC(O,A)(B7,w7) \tkzGetPoints{C7}{D7}   
  \tkzFillCircle[color=white](B7,w7)  
  \tkzDrawArc(B7,C7)(D7)    

  \tkzDefPoint(cos(3.55), sin(3.55)){z8}
  \tkzDefPoint(cos(3.6), sin(3.6)){w8} 
  \tkzDefCircle[orthogonal through=z8 and w8](O,A) \tkzGetPoint{B8}
  \tkzInterCC(O,A)(B8,w8) \tkzGetPoints{C8}{D8}   
  \tkzFillCircle[color=white](B8,w8)  
  \tkzDrawArc(B8,C8)(D8)    

  \tkzDefPoint(cos(3.65), sin(3.65)){z9}
  \tkzDefPoint(cos(3.75), sin(3.75)){w9} 
  \tkzDefCircle[orthogonal through=z9 and w9](O,A) \tkzGetPoint{B9}
  \tkzInterCC(O,A)(B9,w9) \tkzGetPoints{C9}{D9}   
  \tkzFillCircle[color=white](B9,w9)  
  \tkzDrawArc(B9,C9)(D9)   

  \tkzDefPoint(cos(4.5), sin(4.5)){z10}
  \tkzDefPoint(cos(5), sin(5)){w10} 
  \tkzDefCircle[orthogonal through=z10 and w10](O,A) \tkzGetPoint{B10}
  \tkzInterCC(O,A)(B10,w10) \tkzGetPoints{C10}{D10}   
  \tkzFillCircle[color=white](B10,w10)  
  \tkzDrawArc(B10,C10)(D10)

  \tkzDefPoint(cos(5.1), sin(5.1)){z11}
  \tkzDefPoint(cos(5.4), sin(5.4)){w11} 
  \tkzDefCircle[orthogonal through=z11 and w11](O,A) \tkzGetPoint{B11}
  \tkzInterCC(O,A)(B11,w11) \tkzGetPoints{C11}{D11}   
  \tkzFillCircle[color=white](B11,w11)  
  \tkzDrawArc(B11,C11)(D11)

  \tkzDrawCircle[color=black, thin, dotted](O,A)

  \node at (-0.5, 0.7) {$T_\ell$};
  \node at (-0.65, 0.85) {$\ell$};
  

\end{tikzpicture}